\documentclass[3p,times]{elsarticle}

\usepackage[utf8]{inputenc} 

\usepackage{amsmath,amsfonts,amssymb}
\usepackage{geometry} 
\usepackage{graphicx} 

\usepackage{booktabs} 
\usepackage{array} 
\usepackage{paralist} 
\usepackage{verbatim} 
\usepackage{subfig} 

\usepackage{fancyhdr} 
\usepackage{sectsty}
\allsectionsfont{\sffamily\mdseries\upshape} 

\usepackage[nottoc,notlof,notlot]{tocbibind} 
\usepackage[titles,subfigure]{tocloft} 

\usepackage{graphicx}
\usepackage{epstopdf}
\newcommand{\field}[1]{\mathbb{#1}}
\newcommand {\R}{\field{R} }

\newcommand {\eps}{\varepsilon }

\newtheorem{prop}{Proposition}[section]

\newtheorem{rem}[prop]{Remark}

\newtheorem{cor}[prop]{Corollary}

\usepackage{color}
\usepackage{colordvi}
\definecolor{magenta}{rgb}{.5,0,.5}
\definecolor{black}{rgb}{1.0,1.0,1.0}
\definecolor{magenta}{rgb}{.1,0,.3}
\definecolor{gruen}{rgb}{0.2,0.5,.5}
\definecolor{light}{rgb}{ 0.992, 0.961,  0.902}
\definecolor{Tan}{rgb}{ 0.992, 0.9,  0.902}

\newcommand{\komment}[1]{{}}

\usepackage{lineno}

\begin{document}

\begin{frontmatter}

\title{The effect of delay on contact tracing}

\author[TUMma,Helm]{Johannes M\"uller\corref{cor1}}
\ead{johannes.mueller@mytum.de}
\address[TUMma]{
TU M\"unchen,  
Centre for Mathematical Sciences, 
Boltzmannstr.\ 3, 
D-85747 Garching, Germany
}
\address[Helm]{
Helmholtz Center Munich, 
Institute for Computational Biology, 
Ingolst\"adter Landstr. 1, 
D-85764 Neuherberg, Germany
}

\author[TUMma]{Bendix Koopmann}

\cortext[cor1]{Corresponding author}

\begin{abstract}
We consider a model for an infectious disease in the onset of an
 outbreak. We introduce contact tracing incorporating a tracing delay. The effect of randomness in the delay and the 
effect of the length of this delay in comparison to the 
infectious period of the disease respectively to 
a latency period on the effect of tracing, given e.g.\ by the change of the reproduction number, is analyzed. We focus particularly on the effect of randomness in the tracing delay.
\end{abstract}
\end{frontmatter}


\section{Introduction}

Contact tracing resp.\ partner notification programs 
are believed to be of central importance for the control
of many infectious diseases: infected persons are questioned for 
recent potentially infectious contacts. In this way, further 
infected and infectious persons are identified in a targeted way, 
often quite early after infection. These persons can be treated 
and isolated, and the number of further infections 
can be reduced. 
For some emerging infections, data analysis
 indicates that 
contact tracing has proven to be a valuable measure -- e.g.\ in the case of SARS~\cite{klinkenberg2006} 
or Ebola~\cite{rivers2014}. 
For other infections, such as tuberculosis, contacts may be this casual that they are hardly recognized; in these cases it is under debate if contact tracing 
pays~\cite{noah1997,begun2013,fox2013}. 
Still, our understanding  of the effectiveness of contact tracing is incomplete. 
In particular, 
the consequences of the different 
time scales involved -- latent period, 
typical time between contacts, and the delay in the tracing process -- remain unclear.
\par\medskip

As contact tracing depends on the detailed contact structure, it is -- 
in contrast to e.g.\ 
mass screening~-- not immediately clear how to model this method appropriately. Local interactions and
correlations have to be taken into account. 
In recent years, basically two different modelling approaches have been successfully developed. 
The first approach~\cite{eames2002,keeling1999,house2010} relies on a fixed contact graph. 
The infection (as well as contact tracing) spreads via the 
edges of this graph, and is modeled as a stochastic
contact process. Pair approximation 
yields a model consisting of ordinary 
differential equations  (ODEs) that gasps the most important features of the dynamics. The mean value of the stochastic process is more or less 
met by these ODEs. This modeling approach gives in particular good results if the degree of nodes is large. 

The second approach is based on a branching process, 
and in particular used to describe  the onset of 
an outbreak~\cite{mueller2000,mueller2003,mueller2007}. 
On the tree of infecteds (the nodes
are infected individuals, a directed edge points from 
infector to infectee) the tracing process takes place. 
If an individual is discovered, adjacent edges have (independently) a certain probability -- the tracing probability -- to be detected. 
As the underlying graph is directed, it is suggestive 
to define forward tracing (if the infector is discovered,
infectees are traced) and backward tracing (if an  infectee is discovered, the infector is traced). Even 
in the very early papers~\cite{hethcote:gonorrhea} 
this concept
has been developed, 
and it has been discussed if forward- or backward tracing 
is more important.

In addition to these two mathematical approaches, 
a lot of work has been done based on simulation 
 models~\cite{kiss2006,kiss2008} and/or to understand the effect 
of contact tracing for certain diseases like influenza, SARS, tuberculosis or Ebola~\cite{eichner2003,klinkenberg2006,fox2013,begun2013,rivers2014}. \par\medskip

In the present work, we take up the discussion how 
a tracing delay -- the time elapsing between the discovery of an infected individual and the identification of
his/her infector and infectees -- influences the efficiency of contact tracing. Fraser~\cite{fraser2004} and Kiss~\cite{kiss2006} already discussed 
the importance of a latent period for contact tracing: a latent period allows one to detect cases 
before they start to spread the diseases and in this 
makes contact tracing more effective. A tracing delay has the converse effect; persons may spread the infection also during the time 
that elapses between  detection of an index case and their own detection by contact tracing. 
Only few models address this delay explicitly. Klinkenberg et al.~\cite{klinkenberg2006} extends the work
of Fraser et al.~\cite{fraser2004} by a tracing delay. 
Approximations of the next generation operator 
for contact tracing were developed.  
Another approach was used by Shaban et al.~\cite{shaban2008}. 
In that
paper, a fixed contact network is considered (as in 
most pair approximation models), but the authors
focus on the onset of an outbreak and use a branching 
process approximation of the process. 
They only take into
account forward tracing. In principle, their model 
allows for general distributions for latency period and 
tracing delay, but the authors concentrate on the 
special case of exponential distributions. 
Ball et al.~\cite{ball2011,ball2015} take up this idea. They also consider only forward tracing 
but assume a homogeneously mixing population.  The authors formulate a multitype-branching process for detected individuals. 
This approach is mathematically particularly appealing, as the theory 
of branching processes can be used to derive analytical 
results.\par\medskip

In the present work we extend the methods developed
in~\cite{mueller2000} to analyse 
delayed contact tracing  
with forward- and backward tracing. We do this analysis 
first for an epidemic with constant contact- and recovery rates, but extend the ideas also to non-constant rates, 
opening the possibility to also consider an infection with 
latent period. The central technique relies on the 
derivation of the probability 
that an individual is still infectious at a 
given age of infection.  
However, in general it is only possible to solve these equations 
numerically. Approximate solutions are derived for 
small tracing probabilities; also an approximation 
for the reproduction number is given. 
The influence of the timing (latent period and 
tracing delay) is discussed. We are particularly interested
in the question how randomness in the tracing delay affects 
the efficiency of contact tracing as given by the reduction of the reproduction number, if forward- or backward tracing is more important, 
and how the interplay between the time scales involved (tracing delay, 
mean infectious period, and 
latency period)  
influences contact tracing.

\section{Model and Analysis}

We consider a randomly mixing, homogeneous population.  
Note that models assuming a homogeneous 
population may behave differently compared with models that assume an underlying contact graph, in particular of the contact graph is sparse. It depends on the disease 
which approach is more appropriate. 
In order to model contact tracing, we start off with 
an SIS/SIR- type of model, 
and focus on the onset of an epidemic (therefore we do not need to specify if a recovered person will be susceptible again or immune). 
In the long run, SIR and SIS models will, of course, behave differently. 
The contact 
rate is denoted by $\beta$. 
Infected persons recover at rate $\gamma$. With probability $p_{\text{obs}}$ recovered persons 
become index cases and trigger (at recovery) a tracing event. 
Tracing, however, does not take place immediately but with a random 
delay $T$, distributed with density $\phi$. We allow generalized functions for this density,
such that a fixed delay is covered by the model. 
For each contact, the delay is an independent realization; 
it is not the case that the tracing delay 
e.g.\ only depends on the index case.  
Infector and infectees of the index case  have, if they are still infectious 
at the time at which contact tracing actually takes place, 
a probability $p$ to be diagnosed. We consider two different modes: either a traced individual is again an index case 
(recursive tracing), or the contact tracing stops (one-step-tracing).\par\medskip

For technical reasons, we introduce the rates $\alpha=(1-p_{\text{obs}})\gamma$ and $\sigma=p_{\text{obs}}\gamma$, and consider $\alpha$ as the 
spontaneous recovery rate (no diagnosis), and $\sigma$ the recovery rate with direct (not via tracing) 
diagnosis of the infection. Note that we focus on the
onset of the epidemic. Therefore, all contacts of an infected 
individual connect to susceptible individuals. The infection process without tracing is well approximated by 
a linear branching process (along the line of the argument of Ball and Donnelly~\cite{ball:donnelly}, see also~\cite{mueller2000,ball2011}). Only
tracing introduces dependencies between individuals. 
In order to analyze this 
process, we first look at backward tracing, then at forward tracing, and at the end we combine both processes to full tracing.

\subsection{Backward tracing -- recursive mode}
Let us assume that only backward tracing takes place, and no forward tracing. We  furthermore consider 
the recursive
tracing mode.
As usual, convolution of two functions $f$ and $g$ is defined by
$(f\ast g)(a) = \int_0^a f(a-\tau) g(\tau)\, d\tau$. Furthermore, we define
$$f^{\#}(a) = \int_0^af(\tau)\,d\tau = (1\ast f)(t).$$

\begin{prop}\label{backPropConst}
Let $\phi(t)$ denote the distribution of the tracing delay $T$, 
$\kappa_-(a)$ the probability to be infectious after time of infection $a$. Then, 
\begin{eqnarray}
\kappa_-'(a) =
-\kappa_-(a) \left\{
\alpha+\sigma 
+
p\,\beta\,\left[ (\phi\ast (1-\kappa_-))(a)
-\alpha (\phi\ast \kappa_-^{\#})(a)\right]
\right\},\quad\kappa_-(0)=1.\label{exctback}
\end{eqnarray}
\end{prop}
{\bf Proof: } We start off with the relation 
$$\kappa_-'(a) = -\kappa_-(a)\left\{\alpha+\sigma+\mbox{rate of tracing}\,(a)\right\}.$$
In order to obtain the rate of (direct or indirect) 
detection of an infected individual 
with age since infection $a$ , we subtract 
from the total removal rate (or hazard) $-\kappa_-'(a)/\kappa_-(a)$  
the rate of spontaneous removal $\alpha$, 
$$ \frac{-\kappa'_-(a)}{\kappa_-(a)}-\alpha.$$
In order to compute the contribution of 
backward contact tracing to the 
removal rate, we consider the infectees 
(generated at rate $\beta$) that are still infectious 
after  $c$ time units (probability $\kappa_-(c)$)
and are detected at age of the infector~$\tau$, 
$$ \int_0^\tau\beta 
       \left(\frac{-\kappa'_-(c)}{\kappa_-(c)}-\alpha\right) \kappa_-(c)\, dc.
       $$
These individuals increase the removal rate of the infector at age of infection $a$, if the tracing delay 
is precisely $a-\tau$ (probability density $\phi(a-\tau)$), 
and the infector is traced, indeed (probability $p$). That is, 
\begin{eqnarray*}
&&\mbox{rate of  tracing}\,(a)
= p\,\int_0^a\phi(a-\tau)\int_0^\tau\beta 
       \left(\frac{-\kappa'_-(c)}{\kappa_-(c)}-\alpha\right) \kappa_-(c)\, dc \, d\tau\\
&=& p\,\beta \int_0^a\phi(a-\tau)\left((1-\kappa_-(\tau) - \alpha \kappa^{\#}_-(\tau)\right)\, d\tau
= p\,\beta\,\left[ (\phi\ast 1)(a)-(\phi\ast \kappa_-)(a)
-\alpha (\phi\ast \kappa_-^{\#})(a)\right].
\end{eqnarray*}
We obtain the integro-differential equation stated above.\par\qed\par\medskip

Let from now on
$\widehat\kappa(a) = e^{-(\alpha+\sigma) a}$ 
for $a\geq 0$, and $\widehat\kappa(a)=0$ for $a<0$. 
\par\medskip

\begin{prop}
The first order approximation of $\kappa_-(a)$ in $p$ reads 
\begin{eqnarray}
\kappa_-(a) &=& \widehat\kappa(a) - p \,p_{\text{obs}} \,
\beta \, \widehat\kappa(a)\,\,(1\ast\phi\ast(1-\widehat\kappa))(a)
+{\cal O}(p^2).\label{approxback}
\end{eqnarray}
\end{prop}
{\bf Proof: }
We go for a first order approximation. Note that 
$\kappa_-(a)$ does not only depend on $a$ but also 
 on $p$ (and some other parameters that we keep constant). 
For a given $a$, we expand $\kappa_-(a)$ as a power series in $p$, viz.
$$ \kappa_-(a) = \sum_{i=0}^\infty p^i\,\kappa_{-,i}(a).$$
In this formula, the functions $\kappa_{-,i}(a)$ 
do not depend on $p$ any more. 
We replace $\kappa_-(a)$ by this expansion in the integro-differential
equation (\ref{exctback}), 
\begin{eqnarray*}
\sum_{i=0}^\infty p^i\,\kappa_{-,i}'(a) =
-\left(\sum_{i=0}^\infty p^i\,\kappa_{-,i}(a)\right) \left\{
\alpha+\sigma 
+
p\,\beta\,\left[ (\phi\ast (1-\sum_{i=0}^\infty p^i\,\kappa_{-,i}))(a)
-\alpha (\sum_{i=0}^\infty p^i\,\phi\ast \kappa_{-,i}^{\#})(a)\right]
\right\},\quad\sum_{i=0}^\infty p^i\,\kappa_{-,i}(0)=1
\end{eqnarray*}
and equate powers of $p$. We find 
for $p^0$ and $p^1$ 
\begin{eqnarray*}
\kappa_{-,0}'(a) &=& - \kappa_{-,0}(a)(\alpha+\sigma),\qquad \kappa_{-,0}(0) = 1,\\
\kappa_{-,1}'(a) &=& 
 - (\alpha+\sigma)\kappa_{-,1}(a) -
\beta\,\kappa_{-,0}(a)\,\left[ (\phi\ast 1)(a)-(\phi\ast \kappa_{-,0})(a)
-\alpha (\phi\ast \kappa_{-,0}^{\#})(a)\right],\quad\kappa_{-,1}(0)=0.
\end{eqnarray*}
Therefore,
 $\kappa_{-,0}(a) = e^{-(\alpha+\sigma)\, a}=\widehat\kappa(a)$ 
and, with 
$\kappa_{-,0}^{\#}(a) = (1-\widehat\kappa(a))/(\alpha+\sigma)$ 
we obtain that
\begin{eqnarray*}
&& \bigg[(\phi\ast 1)(a)-(\phi\ast \kappa_{-,0})(a)
-\alpha (\phi\ast \kappa_{-,0}^{\#})(a)\bigg]=
(\phi\ast (1-\widehat\kappa))(a)
-\frac\alpha{\alpha+\sigma} (\phi\ast (1-\widehat\kappa))(a)\\
&=&
\frac\sigma{\alpha+\sigma} (\phi\ast (1-\widehat\kappa))(a). 
\end{eqnarray*}
Hence, the first order correction of $\widehat\kappa$ by backward tracing
reads
\begin{eqnarray*}
\kappa_{-,1}(a) &=&-\,\int_0^ae^{-(\alpha+\sigma)(a-\tau)}\,
\frac{\beta\sigma}{\alpha+\sigma} e^{-(\alpha+\sigma)\tau}\,(\phi\ast (1-\widehat\kappa))(\tau)\,d\tau\\
&=&  -\,\frac{\beta\sigma}{\alpha+\sigma} \widehat\kappa(a)\,\,(1\ast\phi\ast(1-\widehat\kappa))(a).
\end{eqnarray*}
Therewith the result follows.\par\qed\par\medskip

\subsubsection{Special case: fixed delay}
Assume that the tracing delay is a deterministic, 
fixed time period $T$. That is, $\phi(a)=\delta_T(a)$. Then, 
$$(\phi\ast (1-\widehat\kappa))(a) = (1-\widehat\kappa(a-T))$$
for $a>T$ and zero else. Hence, for $a>T$, 
\begin{eqnarray*}
(1\ast\phi\ast(1-\widehat\kappa))(a) &=& 
\int_0^a \phi\ast (1-\widehat\kappa)(\tau) \, d\tau
=  \int_T^a (1-\widehat\kappa(a-T)) \, d\tau\\
&=& (a-T) - \frac 1 {\alpha+\sigma}\,\left(1-e^{-(\alpha+\sigma)(a-T)}\right).
\end{eqnarray*}
All in all, we obtain for $a<T$ that 
$$ \kappa_-(a) = e^{-(\alpha+\sigma)a} = \widehat\kappa(a)$$
and for $a>T$
$$ \kappa_-(a) =  \widehat\kappa(a)
-  p \,p_{\text{obs}} \,
\beta  \,\,   \widehat\kappa(a) \,\,\left\{
(a-T) 
- \frac {1- \widehat\kappa(a-T)} {\alpha+\sigma}\,
\right\}+{\cal O}(p^2)
.$$

\begin{figure}
\begin{center}
\includegraphics[width=10cm]{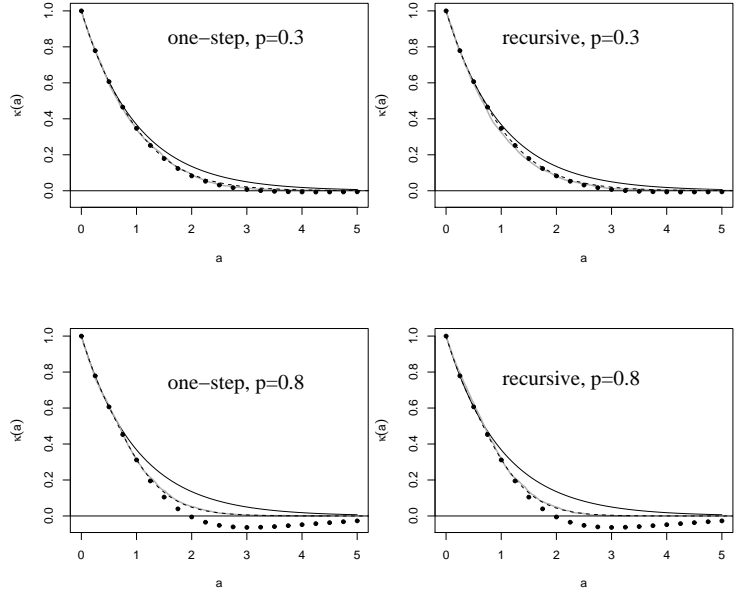}
\end{center}
\caption{$\kappa_-(a)$ for backward tracing (one step, recursive tracing, $p=0.3$, and $p=0.8$, as indicated). Solid line: $\widehat\kappa(a)$, 
gray line: simulated data, dashed line: $\kappa_-(a)$ according to the theory, 
dots: first order approximation of $\kappa_-(a)$. 
Parameters: $\beta=2$, $\alpha=0.1$, $\sigma=0.9$, fixed tracing delay, $T=0.5$.}
\label{figBack}
\end{figure}

The numerical analysis shows an excellent agreement of theory as given by eq.\ (\ref{exctback}) and  
simulation (figure~\ref{figBack}). However, the first order approximation (\ref{approxback}) is
not valid any more if $p$ becomes large.

\subsubsection{Special case: exponential delay}
Apart of a fixed delay, an exponentially distributed 
delay is another natural choice. Let the mean value 
be $T$, 
$\phi(a) = e^{- a/T}/T$. We assume that 
$T(\alpha+\sigma)\not = 1$. Straight forward computations yield 
$$
\kappa_-(a)= \widehat\kappa(a)- p \,p_{\text{obs}} \,
\beta \,  \widehat\kappa(a)\,\left\{
a+ T\,\,\,
\frac{T(\alpha+\sigma)}{1-T(\sigma+\alpha)}\left(
(1-e^{-a/T}) -\frac{1- \widehat\kappa(a)}{(T(\sigma + \alpha))^2}\right)
\right\}+ \mathcal{O}(p^2).
$$
Figure~\ref{effKappaFig} indicates that the exponentially distributed delay has a larger effect than the fixed delay. 
We come later back to this observation and discuss the 
presumable mechanism behind this finding.

\subsubsection{Rates depending on age since infection}
We generalize the model assumptions and allow for 
the case that 
$\beta$, $\alpha$ and $\sigma$ depend on the age
 since infection, e.g., $\beta=\beta(a)$. These relaxed
 assumptions allow to consider the interplay between contact tracing and a latency period.  We only look at the fully recursive case, since the one-step tracing case 
is more simple.\\ 
 The argument here parallels that 
of proposition~\ref{backPropConst}. 
The rate at which an infectee is detected 
at age of infector $\tau$ is now given by 
$$ \int_0^\tau\beta(\tau-c) 
       \left(\frac{-\kappa'_-(c)}{\kappa_-(c)}-\alpha(c)\right) \kappa_-(c)\, dc.
       $$
Hence, the removal rate due to contact tracing 
reads
\begin{eqnarray*}
&&\mbox{rate of removal by contact tracing}(a)\\
&=& p\,\int_0^a\phi(a-\tau)\int_0^\tau\beta(\tau-c) 
       \left(\frac{-\kappa'_-(c)}{\kappa_-(c)}-\alpha(c)\right) \kappa_-(c)\, dc \, d\tau\\
&=& - p\,\int_0^a\phi(a-\tau)\int_0^\tau\bigg(\beta(\tau-c) 
       \kappa'_-(c) +\beta(\tau-c) \alpha(c)\kappa_-(c)\bigg)\, dc \, d\tau.
\end{eqnarray*}
Unfortunately, in general this expression cannot be 
simplified. We obtain the integro-differential equation for $\kappa_-(a)$ 
\begin{eqnarray}
\kappa_-'(a) = 
-\kappa_-(a) \left(
\sigma(a) + \alpha(a)
- p\,\int_0^a\phi(a-\tau)\int_0^\tau\bigg(\beta(\tau-c) 
       \kappa'_-(c) +\beta(\tau-c) \alpha(c)\kappa_-(c)\bigg)\, dc \, d\tau
\right).
\end{eqnarray}

\subsubsection{Backward tracing -- one-step tracing}
We turn to one-step backward tracing. The 
basic argument stays the same as above, but the equations become 
slightly more simple.\par\medskip

\begin{prop}
Let $\phi(t)$ denote the distribution of the tracing delay, 
$\kappa_-(a)$ the probability to be infectious after time of infection $a$, and $\kappa_-^{\#}(a) = \int_0^a\kappa_-(\tau)\,d\tau$. Then, 
\begin{eqnarray}
\kappa_-'(a) =
-\kappa_-(a) \left\{
\alpha+\sigma 
+
p\,\beta\,\sigma\,(\phi\ast \kappa_-^{\#})(a)
\right\},\quad\kappa_-(0)=1.\label{backExaktOneStep}
\end{eqnarray}
\end{prop}
{\bf Proof: } As before, 
$$\kappa_-'(a) = -\kappa_-(a)\left\{\alpha+\sigma+\mbox{rate of tracing}\, (a)\right\}.$$
The rate of (direct) detection as infected individual is just $\sigma$,
that is, 
\begin{eqnarray*}
&&\mbox{rate of  tracing}(a)
= p\,\int_0^a\phi(a-\tau)\int_0^\tau\beta 
       \sigma \kappa_-(c)\, dc \, d\tau
= p\,\beta \,\sigma \int_0^a\phi(a-\tau)\, \kappa^{\#}_-(\tau)\, d\tau
\end{eqnarray*}
We obtain the integro-integral equation stated above.\par\qed\par\medskip

\begin{prop}
The first order approximation in $p$ reads 
\begin{eqnarray}
\kappa_-(a) &=& \widehat\kappa(a) - p \,p_{\text{obs}} \,
\beta \, \widehat\kappa(a)\,\,(1\ast\phi\ast(1-\widehat\kappa))(a)
+{\cal O}(p^2).
\end{eqnarray}
\end{prop}
The proof parallels that of proposition~\ref{approxback}. 
Note that a path of length $2$ has a probability ${\cal O}(p^2)$ to be traced. Hence, the first order approximation of the full recursive backward tracing only takes into account the tracing of immediate neighbors, similar to one-step tracing. This heuristics already indicates that the first order approximation for recursive- and one-step-tracing coincide.\par\medskip

In order to obtain a numerical indication of the quality of our approximation, 
we use $\phi(a) = \delta_T(a)$ as tracing delay. We do not find 
a strong difference between one-step and recursive backward tracing (see figure~\ref{figBack}). As in recursive backward tracing, the first order approximation is well suited for the complete process if the tracing probability is small ($p=0.3$), while for larger tracing probabilities ($p=0.8$) the discrepancy between approximation and exact solution (resp.\ simulations) becomes
more serious.

\subsection{Forward tracing}

Now we proceed to forward tracing. We only discuss recursive forward
tracing in detail, as one-step forward tracing is very similar; 
we will note how to handle one-step forward tracing in remark~\ref{oneStepForward}. Before we formulate the central proposition of
this section, we introduce some 
more notation.\par\medskip

{\bf Definition: }{\it 
Let $\kappa_i(a|b)$ denote the probability for an individual
of generation $i$ 
to be still infectious at age of infection $a$ if the infector has
age of infection $a+b$.}
\par\medskip

\begin{prop}\label{forwardProp}
 We find for $\kappa_i^+(a)$ the recursion formula 
\begin{eqnarray*}
\kappa_{i-1}^+(b)\,\kappa_i(a|b) & = & 
\widehat \kappa(a)\bigg\{
\kappa_{i-1}^+(b)
-p\, \int_0^a 
\left(
{-\kappa_{i-1}^+}'(b+c)
-\,\alpha\,\kappa_{i-1}^+(b+c)\,
\right)\,\int_c^{a}\phi(a'-c)\, da'\, dc\bigg\}\\
\kappa_i^+(a) & = & \frac{\int_0^\infty \kappa_{i}^+(a|b)\kappa_{i-1}^+(b)\, db}{\int_0^\infty \kappa_{i-1}^+(\tau)\,d\tau }.
\end{eqnarray*}
\end{prop}
{\bf Proof: } If the individual has not been traced so far, then its probability to be infectious is that of the zero'th generation, 
$\widehat \kappa(a)$. This probability is decreased by tracing via the infector. 
Hence, to obtain $\kappa^+_i(a|b)$, we multiply $\widehat \kappa(a)$ by 
the probability {\it not} to be traced via the infector. 
This probability is one minus the probability to be
traced. This, in turn, is $p$ times the
probability that during the interval under consideration
a (delayed) tracing event did take place.
\\
The probability for the infector to be still infectious at age of infection $a+b$ is given by
$$\kappa_{i-1}^+(a+b)/\kappa_{i-1}^+(a)$$
as we know at age $a$ an infectious event had happened. 
The rate at which the infector is observed at age $b+c$ is given by
$$
\frac{\kappa_{i-1}^+(b+c)}{\kappa_{i-1}^+(b)}\,\left(
\frac{-{\kappa_{i-1}^+}'(b+c)}{\kappa_{i-1}^+(b+c)}-\alpha
\right)
=
\left(
\frac{-{\kappa_{i-1}^+}'(b+c)}{\kappa_{i-1}^+(b)}
-\,\frac{\alpha\,\kappa_{i-1}^+(b+c)}{\kappa_{i-1}^+(b)}\,
\right)
$$
Therefore, the rate at which 
the infector triggers a tracing event at 
$a'\in[0,a)$ reads 
$$\int_0^{a'} 
\left(
\frac{-{\kappa_{i-1}^+}'(b+c)}{\kappa_{i-1}^+(b)}
-\,\frac{\alpha\,\kappa_{i-1}^+(b+c)}{\kappa_{i-1}^+(b)}\,
\right)\phi(a'-c)\, dc$$
Since we only want to know if an tracing event has been triggered 
before age $a$ of infection 
we integrate over $a'$, and find
\begin{eqnarray*}
&&\int_0^a\,\int_0^{a'} 
\left(
\frac{-{\kappa_{i-1}^+}'(b+c)}{\kappa_{i-1}^+(b)}
-\,\frac{\alpha\,\kappa_{i-1}^+(b+c)}{\kappa_{i-1}^+(b)}\,
\right)\phi(a'-c)\, dc\, da'\\
&=&
\int_0^a 
\left(
\frac{-{\kappa_{i-1}^+}'(b+c)}{\kappa_{i-1}^+(b)}
-\,\frac{\alpha\,\kappa_{i-1}^+(b+c)}{\kappa_{i-1}^+(b)}\,
\right)\,\int_c^{a}\phi(a'-c)\, da'\, dc.
\end{eqnarray*}
Therefore,
$$\kappa_i^+(a|b) 
= 
\widehat \kappa(a)\,\,
\left\{
1 - p 
\int_0^a 
\left(
\frac{-{\kappa_{i-1}^+}'(b+c)}{\kappa_{i-1}^+(b)}
-\,\frac{\alpha\,\kappa_{i-1}^+(b+c)}{\kappa_{i-1}^+(b)}\,
\right)\,\int_c^{a}\phi(a'-c)\, da'\, dc
\right\}. 
$$
If we multiply this equation by $\kappa_{i-1}^+(b)$, we obtain the first
equation of our proposition. 
As the distribution of the age-since-infection of the infector at an
infectious event is given by
$$ \frac{\kappa_{i-1}^+(a)}{\int_0^\infty \kappa_{i-1}^+(\tau)\,d\tau }$$
also the second equation holds true.
\par\qed\par\medskip

\begin{rem}\label{oneStepForward}
In order to obtain the parallel formula for one step tracing, 
we replace the rate for direct and indirect detection $$\left(
\frac{-{\kappa_{i-1}^+}'(b+c)}{\kappa_{i-1}^+(b+c)}-\alpha
\right)$$ by the rate for direct detection only, $\sigma$, 
and find in this way the recursive formula 
\begin{eqnarray*}
\kappa_{i-1}^+(b)\,\kappa_i(a|b) & = & 
\widehat \kappa(a)\bigg\{
\kappa_{i-1}^+(b)
-p\, \sigma\,\int_0^a 
\kappa_{i-1}^+(b+c)\,
\,\int_c^{a}\phi(a'-c)\, da'\, dc\bigg\}\\
\kappa_i^+(a) & = & \frac{\int_0^\infty \kappa_{i}^+(a|b)\kappa_{i-1}^+(b)\, db}{\int_0^\infty \kappa_{i-1}^+(\tau)\,d\tau }.
\end{eqnarray*}
\end{rem}

Now we return to recursive forward tracing. 
In order to obtain a 
first order approximation of $\kappa_i^+(a)$, 
we note that the recursion formula can be written
as 
\begin{eqnarray*}
\kappa_{i}^+(a)\, & = & 
\widehat\kappa(a)\left\{
1
-p\, \frac{
\int_0^\infty\,\,\int_0^a 
\left(
{-\kappa_{i-1}^+}'(b+c)
-\,\alpha\,\kappa_{i-1}^+(b+c)\,
\right)
\,\int_c^{a}\phi(a'-c)\, da'\, dc
\,db
}{
\int_0^\infty \kappa_{i-1}^+(b)\,db
}
\right\}.
\end{eqnarray*}
For a first order approximation of $\kappa_i^+$, 
only a zero order approximation of $\kappa_{i-1}^+$ is
required. We know 
$\kappa_{i-1}^+(a)=\widehat\kappa(a)+ {\cal O}(p)=\exp(-(\sigma+\alpha)a)+ {\cal O}(p)$.
Accordingly, we find for the 
appropriate approximation of the integral expression 
\begin{eqnarray*}
&&
\int_0^\infty\,\,\int_0^a 
\left(-{\kappa_{i-1}^+}'(b+c)
-\,\alpha\,\kappa_{i-1}^+(b+c)\right)\,
\,\int_c^{a}\phi(a'-c)\, da'\, dc
\,db\\
&=&
\int_0^\infty\,\,\int_0^a 
\,\sigma\,e^{-(\alpha+\sigma)(b+c)}\,
\int_c^{a}\phi(a'-c)\, da'\, dc
\,db + {\cal O}(p)\\
&=&\,\sigma\,\,
\int_0^\infty e^{-(\alpha+\sigma)b}\,db\,\,
\int_0^a 
\,e^{-(\alpha+\sigma)c}\,
\,
\int_c^{a}\phi(a'-c)\, da'\, dc+ {\cal O}(p)\\
&=&
\frac{-\sigma}{(\sigma+\alpha)^2}\,
\,\int_0^a 
\,\frac d {dc}e^{-(\alpha+\sigma)c}\,
\,
\int_0^{a-c}\phi(a'')\, da''\, dc+ {\cal O}(p)\\
&=&\frac{-\sigma}{(\sigma+\alpha)^2}\,
\left\{
e^{-(\alpha+\sigma)c}\,
\,
\int_0^{a-c}\phi(a')\, da'
\bigg|_{c=0}^a
-
\,\int_0^a 
\,e^{-(\alpha+\sigma)c}\,
\frac d {dc}\,
\int_0^{a-c}\phi(a')\, da'\, dc
\right\}+ {\cal O}(p)\\
&=&\frac{-\sigma}{(\sigma+\alpha)^2}\,
\left\{
-
\int_0^{a}\phi(a')\, da'
+
\,\int_0^a 
\,e^{-(\alpha+\sigma)c}\,
\phi(a-c)\,  dc
\right\}+ {\cal O}(p)\\
&=&\frac{-\sigma}{(\sigma+\alpha)^2}\,
\,
\big((1-\widehat\kappa)\ast\phi\big)(a)+ {\cal O}(p)
\end{eqnarray*}
where we assumed that $\int_0^\eps \phi(a')\, da'\rightarrow 0$ for $\eps\rightarrow 0$; this 
may
not the true if a fraction of index cases induce an immediate contact tracing event. 
As $\int_0^\infty \kappa_{i-1}^+(b)\,db = \int_0^\infty \widehat\kappa(b)\,db
+{\cal O}(p)=1/(\alpha+\sigma)+{\cal O}(p)$, and $p_{\text{obs}}=\sigma/(\sigma+\alpha)$, 
we obtain the following corollary.
\par\medskip

{\bf Corollary } {\it 
If 
 $\int_0^\eps \phi(a')\, da'\rightarrow 0$ for $\eps\rightarrow 0$, then 
the first order approximation for 
$\kappa_i(a)$ is independent on $i$ (for $i>0$) 
and reads
$$ \kappa_i(a) = \widehat\kappa(a) - 
\,p \,\, p_{\text{obs}}\,\, \widehat\kappa(a)\,\,
\big((1-\widehat\kappa)\ast\phi\big)(a) 
+{\cal O}(p^2).
$$
}
Note that the first order approximation coincides with $\kappa_1^+(a)$.

\subsubsection{Special case: fixed delay}
\begin{figure}
\begin{center}
\includegraphics[width=10cm]{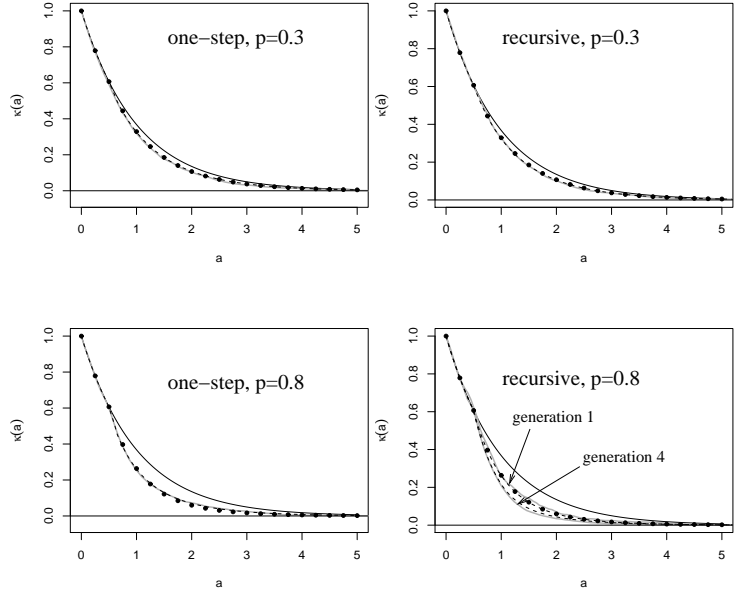}
\end{center}
\vspace*{-0.75cm}
\caption{$\kappa^+_i(a)$ for forward tracing (one step, recursive tracing, $p=0.3$, and $p=0.8$, as indicated). Solid line: $\widehat\kappa(a)$, 
gray line: simulated data for generation 4 (if not indicated differently), dashed line: $\kappa^+_i(a)$ according to the theory  for generation 4 (if not indicated differently),
dots: first order approximation of $\kappa^+_i(a)$. 
Parameters: $\beta=2$, $\alpha=0.1$, $\sigma=0.9$, fixed tracing delay, $T=0.5$.}
\label{figForward}
\end{figure}

If we have a fixed delay, that is $\phi(a)=\delta_T(a)$, then
$$\phi\ast (1-\widehat\kappa)(a)= (1-\widehat\kappa(a-T))$$
for $a>T$ and zero else. Hence, for $a>T$, 
$$ 
 \kappa_i(a) = \widehat \kappa(a) - 
\,\,p \,\, p_{\text{obs}}\,\, \widehat\kappa(a)\,\,
(1-\widehat\kappa(a-T))\,
+{\cal O}(p^2),
$$
 and $ \kappa_i(a) = \widehat\kappa(a) +{\cal O}(p^2)$ else.
Numerical simulations indicate that the first order approximation 
is well suited for the stochastic process, even if the tracing 
probability becomes larger (see figure~\ref{figForward}).

\subsubsection{Special case: exponential delay}
\begin{figure}
\begin{center}
\includegraphics[width=12cm]{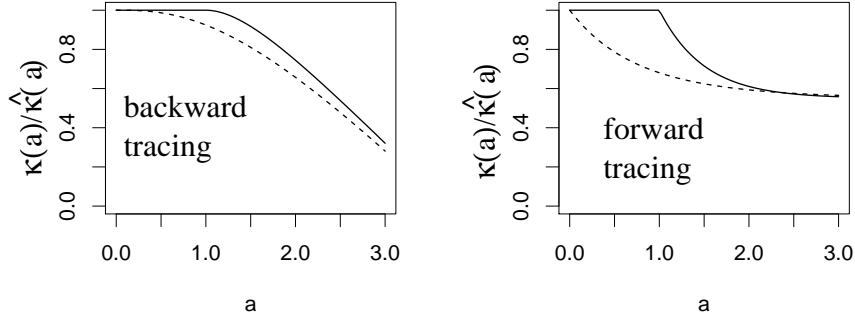}
\end{center}
\caption{First order approximation of $\kappa^\pm(a)$ 
over $\widehat\kappa(a)$; left panel: backward-tracing, 
right panel: forward-tracing. Solid line: $\theta(a)=\delta_T(a)$, 
dashed line: $\phi(a)=e^{-a/T}/T$. 
Parameters: 
$\sigma = \alpha = 1$, 
$p     = 0.3$, 
$\beta = 3$, and 
$T = 1$.}\label{effKappaFig}
\end{figure}

In order to compare fixed and random delay, we 
choose for $\phi(a)$ an exponential distribution with expectation $T$, 
$\phi(a)=e^{-a/T}/T$and obtain
$$\kappa_i(a)=\widehat{\kappa}(a)-\,p\, p_{\text{\text{obs}}}\, \widehat{\kappa}(a)\left(
	1-e^{-a/T} 														-\frac{ \widehat{\kappa}(a) - e^{-a/T}}{1 -T(\alpha +\sigma)}
\right) + \mathcal{O}(p^2).
$$
As before, the exponentially distributed delay 
induces a higher
effect than the fixed delay (Figure~\ref{effKappaFig}).

\subsubsection{Rates depending on age since infection}

We again generalize the considerations above to the
case that $\beta$, $\alpha$, and $\sigma$ depend on 
the time since infection. It is straight to obtain 
the equation for $\kappa_{i}(a)$, but unfortunately,
the equations become even more unhandy as those above. 
However, a first order analysis in $p$ is possible
and yields useful results. 
\par\medskip
The arguments completely parallel that of  proposition~\ref{forwardProp}. We find 
for $\kappa_i^+(a)$ the recursion formula 
\begin{eqnarray*}
\kappa_{i-1}^+(b)\,\kappa_i(a|b) & = & 
\widetilde\kappa(a)\bigg\{
\kappa_{i-1}^+(b)
-p\, \int_0^a 
\left(
{-\kappa_{i-1}^+}'(b+c)
-\,\alpha(b+c)\,\kappa_{i-1}^+(b+c)\,
\right)\,\int_c^{a}\phi(a'-c)\, da'\, dc\bigg\}\\
\kappa_i(a) & = & \frac{\int_0^\infty \kappa_{i}^+(a|b)\beta(b)\kappa_{i-1}^+(b)\, db}{\int_0^\infty \beta(\tau)\kappa_{i-1}^+(\tau)\,d\tau }\\
\widetilde\kappa(a) &=& e^{\int_0^a\,\sigma(\tau)+\alpha(\tau)\, d\tau}.
\end{eqnarray*}

\subsection{Full Tracing}
Let $\kappa_i(a)$ denote the probability to be infectious at age
of infection $a$. Its straight forward to combine forward and backward tracing,
as we only need to repeat the argumentation of the last section, 
but taking into account that $\kappa_0(a)$ is for full tracing not
given by $\kappa_0^+(a) = \exp(-(\alpha+\sigma)a)$, but by $\kappa_-(a)$. Hence we have
the following result.

\begin{prop} We find for $\kappa_0(a)$ the 
\begin{eqnarray}
\kappa_0'(a) =
-\kappa_0(a) \left\{
\alpha+\sigma 
+
p\,\beta\,\left[ (\phi\ast (1-\kappa_0))(a)
-\alpha (\phi\ast \kappa_0^{\#})(a)\right]
\right\},\quad\kappa_0(0)=1.\label{exctfulla}
\end{eqnarray}
and  for $\kappa_i(a)$ for $i>0$ the recursion formula 
\begin{eqnarray}
\kappa_{i-1}(b)\,\kappa_i(a|b) & = & 
\kappa_0(a)\bigg\{
\kappa_{i-1}(b)
-p\, \int_0^a 
\left(
{-\kappa_{i-1}}'(b+c)
-\,\alpha\,\kappa_{i-1}(b+c)\,
\right)\,\int_c^{a}\phi(a'-c)\, da'\, dc\bigg\}\label{exctfullb}\\
\kappa_i(a) & = & \frac{\int_0^\infty \kappa_{i}^+(a|b)\kappa_{i-1}(b)\, db}{\int_0^\infty \kappa_{i-1}(\tau)\,d\tau }.\label{exctfullc}
\end{eqnarray}
\end{prop}
These formulas are exact but not handy. For small
tracing probabilities, however, a first order approximation is enough 
to estimate the effect of contact tracing. 
In order to do so we again only need to put together the results 
obtained for forward- and backward tracing.

\begin{prop}
The first order approximation in $p$ reads 
\begin{eqnarray}
\kappa_0(a) &=& \widehat\kappa(a)
\bigg\{1 - p \,\, p_{\text{obs}}\,\,\beta\,\, 
(1\ast\phi\ast(1-\widehat\kappa))(a)
\bigg\}
+{\cal O}(p^2).
\end{eqnarray}
If 
 $\int_0^\eps \phi(a')\, da'\rightarrow 0$ for $\eps\rightarrow 0$, then 
the first order approximation for 
$\kappa_i(a)$ is independent on $i$ (for $i>0$) 
and reads 
\begin{eqnarray}
\kappa_i(a) 
 &=& \widehat\kappa(a)
\bigg\{1 
- p \,\, p_{\text{obs}}\,\,\beta\,\, 
(1\ast\phi\ast(1-\widehat\kappa))(a)
-
 \,\,p \,\, p_{\text{obs}}\,\, 
\big(\phi\ast(1-\widehat\kappa)\big)(a) \bigg\}
+{\cal O}(p^2).
\end{eqnarray}
\end{prop}

\begin{rem}
The reproduction number 
in the $i$'th generation 
with contact tracing simply reads
\begin{eqnarray*}
&&R^{(i)} =  \int_0^\infty\beta\kappa_i(a)\, da\\
& = & R_0 
 - p \,\, p_{\text{obs}}\,\,
 \left\{\int_0^\infty\beta^2\,\, 
\widehat\kappa(a)\,\,(1\ast\phi\ast(1-\widehat\kappa))(a)\, da
 +
 \int_0^\infty\beta\widehat\kappa(a)\,\,
\big((1-\widehat\kappa)\ast\phi\big)(a) \, da\right\}
+{\cal O}(p^2).
\end{eqnarray*}
In this approximation, the effects of forward- 
and backward tracing are clearly separated.
\end{rem}

\subsubsection{Fixed delay}

\begin{figure}
\begin{center}
\includegraphics[width=10cm]{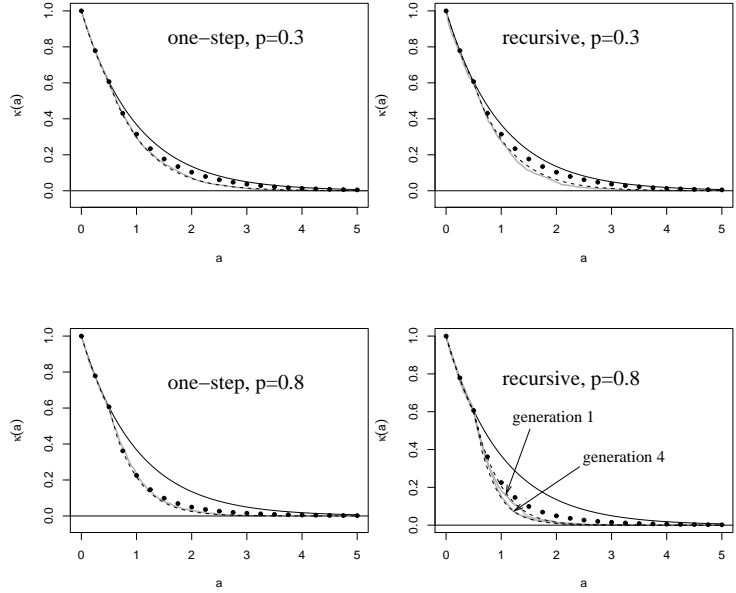}
\end{center}
\vspace*{-0.75cm}
\caption{$\kappa_i(a)$ for full tracing (one step, recursive tracing, $p=0.3$, and $p=0.8$, as indicated). Solid line: $\widehat\kappa(a)$, 
gray line: simulated data for generation 4 (if not indicated differently), dashed line: $\kappa_i(a)$ according to the theory  for generation 4 (if not indicated differently),
dots: first order approximation of $\kappa_i(a)$.}
\label{figFull}
\end{figure}

As indicated, the first order approximation for 
$\kappa_i(a)$ is just a combination of the approximations 
for forward- and backward tracing. We obtain 
for $a<T$ that $\kappa_i(a)=\widehat\kappa(a) + {\cal O}(p^2)$ and for $a>T$
$$ \kappa_-(a) = e^{-(\alpha+\sigma)a}
-  p \,p_{\text{obs}} \,
  \,\,  \widehat \kappa(a) \,\,\left\{
\beta(a-T) 
- \frac \beta {\alpha+\sigma}\,\big(1-\widehat \kappa(a-T)\big)
+ (1-\widehat\kappa(a-T))
\right\}+{\cal O}(p^2)
.$$
As before, numerical simulations indicate that the first order approximation 
is well suited as long as the 
probability is not too large (see figure~\ref{figFull}).

We find 
\begin{eqnarray*}
\int_T^\infty \hat\kappa(a) (1\ast\phi\ast(1-\widehat\kappa))(a)\, da
&=&\widehat\kappa(T)\,
\int_T^\infty \hat\kappa(a-T) \left\{(a-T) - \frac 1 {\alpha+\sigma}\,\left(1-e^{-(\alpha+\sigma)(a-T)}\right)\right\}\, da\\
& = & \frac{\widehat \kappa(T)}{2(\alpha+\sigma)^2}.
\end{eqnarray*}
And, with 
$\phi\ast (1-\widehat\kappa)(a)= (1-\widehat\kappa(a-T))$, 
we conclude that
\begin{eqnarray*}
\int_T^\infty \widehat\kappa(a)\,\, (\phi\ast (1-\widehat\kappa))(a)\, da
&=&
\widehat\kappa(T)  
\int_0^\infty \widehat\kappa(a)  (1-\widehat\kappa(a))\, da = \frac {\widehat\kappa(T)} {2(\alpha+\sigma)}. 
\end{eqnarray*}
Hence, we obtain the following proposition.

\begin{prop} For a fixed delay, $\phi(a)=\delta_T(a)$, we obtain for $i>0$ 
\begin{eqnarray} 
R_{ct} := R^{(i)} = R_0 
- \frac 1 2\,\,p\,p_{\text{obs}}\, \widehat\kappa(T) R_0(R_0+1)+{\cal O}(p^2).
\end{eqnarray}
\end{prop}
The effect of contact tracing is (in first order of $p$) exponentially decreasing in the
tracing delay. The time scale of this exponential decrease is given by the total removal rate $\gamma=\alpha+\sigma$.

\subsubsection{Special case: exponential delay}

We choose $\phi(a)=e^{-a/T}/T$. We do not state the 
first order approximation of $\kappa_i(a)$, as it is only necessary to 
combine
the corresponding results from forward- and backward 
tracing. We focus on the first order effect 
of contact tracing on  
the reproduction number. The effect due to 
backward tracing is described by 
\begin{eqnarray*}
&&\beta \int_0^\infty p \,p_{\text{obs}} \,
\beta \,  \widehat\kappa(a)\,\left\{
a+ T\,\,\,
\frac{T(\alpha+\sigma)}{1-T(\sigma+\alpha)}\left(
(1-e^{-a/T}) -\frac{1- \widehat\kappa(a)}{(T(\sigma + \alpha))^2}\right)
\right\}\, da
=
 \frac{p\,p_{\text{obs}} \, R_0^2}{2(1+T(\alpha+\sigma))}
\end{eqnarray*}
and that for forward tracing
\begin{eqnarray*}
&&\beta\int_0^\infty
\,p\, p_{\text{\text{obs}}}\, \widehat{\kappa}(a)\left(
	1-e^{-a/T} 														-\frac{ \widehat{\kappa}(a) - e^{-a/T}}{1 -T(\alpha +\sigma)}
\right)
\, da
= \frac {p\, p_{\text{\text{obs}}}\,R_0\, }{2(1+T(\alpha+\sigma))}.
\end{eqnarray*}

\begin{prop} For an exponentially distributed delay, 
$\phi(a)=e^{-a/T}/T$, we obtain for $i>0$ 
\begin{eqnarray} 
R_{ct} := R^{(i)} = R_0 
- \frac 1 2\,\,p\,p_{\text{obs}}\, \frac{ R_0(R_0+1)}{1+T(\alpha+\sigma)}+{\cal O}(p^2).
\end{eqnarray}
\end{prop}
We again find that the exponential distributed delay 
yields a higher effect than the fixed delay (if we 
compare distributions with the same expectation). 
While the effect of the fixed delay decrease exponentially in 
$T$, the exponential delay only yields a polynomial decay. 
Indeed, $1/(1+T(\alpha+\sigma))$ is the [0/1]-Pad\'e approximation of 
$\widehat\kappa(T)=\exp(-T(\alpha+\sigma))$. 
It is remarkable that in both 
cases, 
the effect of the delay only depends on $T(\alpha+\sigma)$,
that is, on the quotient of the expected delay over the 
expected time of infection (in absence of contact tracing).

\subsubsection{Rates depending on age since infection}
To obtain the full model for the case if the rates depend on the age since infection, we again only have to combine 
forward- and backward tracing for this case. 
All in all, we obtain the equations for the probability 
to be infective at age of infection $a$ for an infected
individual of the $i$'th generation $\kappa_i(a)$, 
\begin{eqnarray}
\kappa_0'(a) &=& 
-\kappa_0(a) \bigg(
\sigma(a) + \alpha(a)\nonumber\\
&&\qquad \quad
- p\,\int_0^a\phi(a-\tau)\int_0^\tau\bigg(\beta(c) 
       \kappa'_0(c) +\beta(c) \alpha(c)\kappa_0(c)\bigg)\, dc \, d\tau
\bigg)\label{exctfullvarba}\\
\kappa_{i-1}(b)\,\kappa_i(a|b) & = & 
\kappa_{0}(a)\bigg\{
\kappa_{i-1}(b)\nonumber\\
&&\qquad
-p\, \int_0^a 
\left(
{-\kappa_{i-1}}'(b+c)
-\,\alpha(b+c)\,\kappa_{i-1}(b+c)\,
\right)\,\int_c^{a}\phi(a'-c)\, da'\, dc\bigg\}\label{exctfullvarb}\\
\kappa_i(a) & = & \frac{\int_0^\infty 
\kappa_{i}(a|b)\beta(b)\kappa_{i-1}(b)\, db}
 {\int_0^\infty \beta(\tau)\kappa_{i-1}(\tau)\,d\tau }.\label{exctfullvarc}
\end{eqnarray}

\subsubsection{Error analysis}

We did focus on an approximate technique: though we are able to 
derive exact equations for $\kappa_i(a)$, these equations are too 
complex to be solved, and hence we basically focus on a first order
approximation in the tracing probability $p$. However, in comparison 
with ``real world epidemics'' this is not the only simplification. 
The model itself is a simplification (SIS/SIR within an unstructured
population), and also the branching process with tracing is
only valid for the onset of the epidemic. 
Depletion of the class of susceptibles and tracing via contacts between 
infected individuals (no transmission of infection happen 
in these contacts) are neglected. Both effects gain importance 
if the disease approaches an endemic state. \par\medskip

Truncation error: In Figures~1 and~2, we compare simulations of the
stochastic process resp.\ the numerical solution of the exact equations with the approximate solutions. 
For these simulations, the basic reproduction number $R_0=2$, 
and $p_{obs}=0.9$. While a reproduction of two is in an realistic range, the fraction of observed cases $p_{obs}=0.9$ is rather high to 
focus on 
the effect of contact tracing as well as to uncover potential 
approximation errors. 
We find by visual inspection 
that for reasonable tracing probabilities ($p$ below $0.3$, say), 
the agreement of our approximation and simulations is satisfying. Only 
if the tracing probability becomes larger ($p\approx 0.8$), 
the error is noticeable. Here, in particular 
the fraction of cases that become an index case $p_{obs}$ plays a 
role: if $p_{obs}$ is small, then even large tracing probabilities 
do not play a role.  \par\medskip

\begin{figure}
\begin{center}
\includegraphics[width=\textwidth]{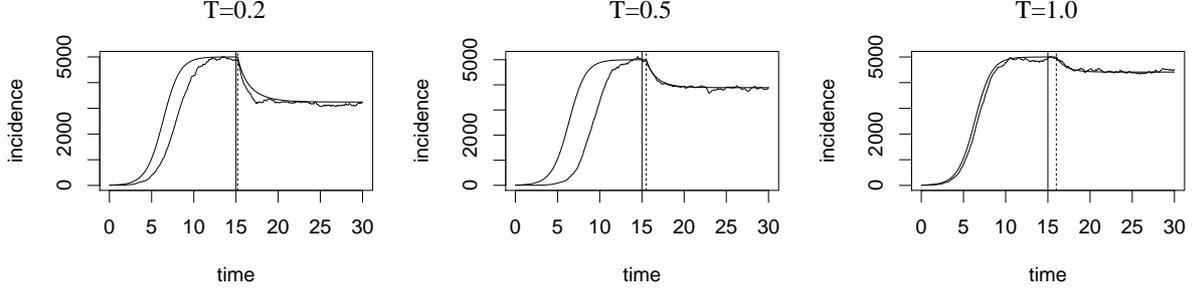}
\end{center}
\vspace*{-0.75cm}
\caption{Simulation of the stochastic (non-smooth line) 
and the approximate deterministic (smooth line)
model. We use a fixed delay $T$, $\beta=2$, $\alpha=0.2$, $\sigma=0.9$,
and recursive tracing.
For $t<15$, $p=0$. At $t=15$ (indicated by the vertical solid line) 
$p$ is increased to $0.3$. The tracing delay is indicated by 
the difference between the dashed and the solid vertical line; at the solid line, contact tracing is introduced, at the 
dashed line, the first time 
tracing effects appear.
}
\label{nonLin}
\end{figure}

Error due to the branching approximation: Reduction of the abundance 
of susceptibles due to the spread of the disease reduces the number of 
secondary cases, and in this, the number of infectious persons that 
can be detected by contact tracing. In order to have at least a heuristic 
method to deal this source of error, we propose to approximate  
what may be called the effective 
removal rate. This removal rate should induce in average the same mean infectious period as the stochastic process. Let us assume that the relative number 
of susceptibles is constant $u$ over a relatively long time period.   
The rate of infectious contacts is reduced from $\beta$ to $\beta \, u$. 
For a fixed delay, the effective reproduction number for a contact 
rate $\beta\, u$ is given by 
$$ R_{ct}\approx \frac{\beta u}{\alpha+\sigma}\left\{1-\frac 1 2 \widehat\kappa(T)(\beta u/(\alpha+\sigma)+1)\right\}.$$
This, in turn, is equivalent with a mean recovery rate $\gamma=\gamma(u)$ given by
$$ \gamma_{eff}(u) = \frac{\beta u}{R_{ct}}
\approx \frac {\alpha+\sigma}{1-\frac 1 2 \widehat\kappa(T)(\beta u/(\alpha+\sigma)+1)}.$$
We neglect in our considerations that also contacts between 
two infected persons take place. Also these contacts 
may lead to tracing events. However, 
if we compare stochastic simulations of the full epidemic process with contact tracing on the one hand, and a deterministic SIS-model with the 
nonlinear recovery rate given above, we find a satisfying agreement (see Fig.~\ref{nonLin}). 
The errors introduced by the saturation of the epidemic process are well
met by the heuristic formula for the effective recovery rate given here.

\subsubsection{The interplay between tracing delay and latency period}

The equation above is too complex to be directly useful in the sense that we obtain deeper inside 
into the interplay of the timing of the disease on the one hand, and tracing on the other hand.
Therefore we concentrate on a special
case: a fixed latency period $T_i$ and a fixed 
tracing delay $T$, 
$$
\beta(a) = \chi_{a>T_i}\,\beta, \quad 
\alpha(a) = \chi_{a>T_i}\,\alpha, \quad 
\sigma(a) = \chi_{a>T_i}\,\sigma,\quad
\phi(a)=\delta_T(a)
$$
where, as usual, $\chi_{a>T_i}$ denotes the characteristic 
function ($\chi_{a>T_i}=1$ if $a>T_i$ and $0$ else). 
We emphasize that $\beta(a)$ denotes a function,
while $\beta$ is a constant (the same for $\alpha$ and $\sigma$). 
 A first
order approximation of the reproduction number is straightforward,
but tedious; we move the computations to~\ref{latenecyAppendix},
and only present the result here. Let $R_{ct}$ denote the reproduction 
number with contact tracing if the number of generations tends to infinity,
and $R_0$ the reproduction number without contact tracing. 
\begin{cor}
\begin{eqnarray}
R_{ct} = 
R_0
- \frac 1 2\,\,p \, p_{\text{obs}}\, R_0\, 
\left[
R_0
\,\,\widehat\kappa(T+T_i)
+
\frac{\widehat\kappa(\max(T_i, T))}
{\widehat\kappa(T_i)}\,\,
\left(
2
-\,
\frac{\widehat\kappa(\max(T_i, T))}
{\widehat\kappa(T)}
\right)
\right]
+{\cal O}(p^2)
\end{eqnarray}
\end{cor}

The tracing effect is monotonously decreasing in the
tracing delay $T$. 
The first order effect consists of two parts, that for backward tracing (see Appendix~\ref{latenecyAppendixBack})
$$ \frac 1 2\,\,p \, p_{\text{obs}}\, R_0^2
\,\,\widehat\kappa(T+T_i)$$
and that for forward tracing (see Appendix~\ref{latenecyAppendixForw})
$$\frac 1 2\,\,p \, p_{\text{obs}}\, R_0\, \,
\frac{\widehat\kappa(\max(T_i, T))}
{\widehat\kappa(T_i)}\,\,
\left(
2
-\,
\frac{\widehat\kappa(\max(T_i, T))}
{\widehat\kappa(T)}
\right).
$$
The backward tracing part is simply exponentially decreasing 
in the tracing delay and the latent period. 

\begin{figure}
\begin{center}
\includegraphics[width=15cm]{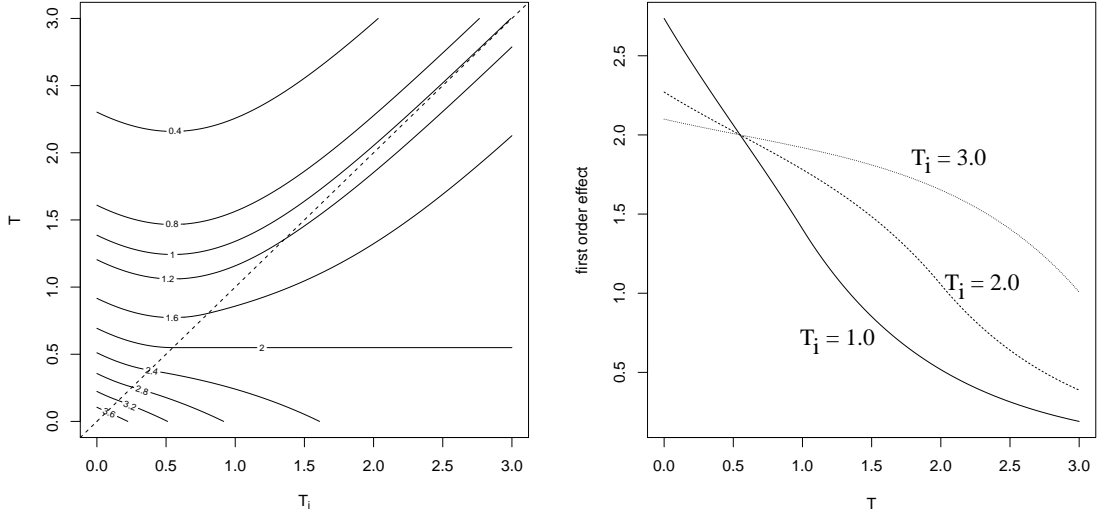}
\end{center}
\caption{Strength of first order effect 
$R_0
\,\,\widehat\kappa(T+T_i)
+
\frac{\widehat\kappa(\max(T_i, T))}
{\widehat\kappa(T_i)}\,\,
\left(
2
-\,
\frac{\widehat\kappa(\max(T_i, T))}
{\widehat\kappa(T)}
\right)$ in dependency on $T$ and $T_i$ (measured in multiples of $1/(\alpha+\sigma)$, the mean infectious period) for $R_0=2$. Left panel: contour plot. Right panel: effect over $T$, $T_i$ fixed at indicated values.}\label{firstOrder}
\end{figure}

The dependency of the froward tracing effect 
on tracing delay resp.\
latent period is more complex. 
First of all, we find that for $T_i=T$, the two delays
in the forward effect cancel each other: There is a race between
tracing and infection. If both mechanisms are subject to 
the same delay, 
the effects of the delays do cancel. If $T<T_i$, the forward tracing
effect is even larger than that without any delay; only if $T>T_i$, 
the effect decreases exponentially in $T$. 
\par\medskip

In order to compare the relative importance of forward- and
backward tracing, we distinguish the cases $T>T_i$ and $T<T_i$.\\
Case $T>T_i$:\\
\begin{eqnarray*}
R_{ct} = 
R_0
- \frac 1 2\,\,p \, p_{\text{obs}}\, R_0\, \widehat\kappa(T+T_i))
\bigg[
R_0
\,\,
+
1\,
\bigg]
+{\cal O}(p^2)
\end{eqnarray*}
As the backward term incorporates $R_0$, and the forward term $1$,
the backward tracing will contribute considerable more to the
(first order) tracing effect.\par\medskip
Case $T<T_i$:\\
\begin{eqnarray*}
R_{ct} = 
R_0 - 
 \frac 1 2\,\,p \, p_{\text{obs}}\, R_0\, 
\widehat\kappa(T+T_i)\,\left[
R_0
\,\,
+
\frac{
2
-\,\widehat\kappa(T_i-T)
}{\widehat\kappa(T+T_i)}
\right]
+{\cal O}(p^2)
\end{eqnarray*}
This time we compare $R_0$ with $(2-\widehat\kappa(T_i-T))/\widehat\kappa(T_i+T)$. In principle, the latter term
can be arbitrarily large, in particular if $T_i$ is large and
$T$ small;  
note that $(2-\widehat\kappa(T_i-T))/\widehat\kappa(T_i+T)$ is always increasing 
in $T_i$ (for $T$ fixed).
For long latency periods and small tracing delays, 
forward tracing gains increasingly  importance. 
 Forward tracing 
benefits from a long latency period. 
\par\medskip

If we fix $T_i$, the first order effect is simply decreasing 
in $T$. If $T_i$ is rather large, a small tracing delay affects
the effect only weakly (see Fig.~\ref{firstOrder}); only if the latency period is 
small or the tracing delay is in the same magnitude (or larger) than
the latency period, the delay in the tracing process strongly 
decreases the effect.\\
If $T$ is fixed, then the effect is in 
general non-monotonously in $T_i$: If $T$ is below  
a certain value, the effect is simply decreasing in $T_i$, if $T$ is above this value, the effect is in a first interval decreasing, but eventually
increasing in $T_i$. The reason for this observation is based on the effect that we observe
the sum of forward- and backward tracing. 
Backward tracing is always decreasing in $T_i$ 
if $T$ is fixed. Forward tracing, however, 
is decreasing if $T_i<T$, but increasing for
$T_i>T$. If the impact of forward tracing is large enough, a non-monotone effect (in $T_i$) may appear.

\section{Discussion}

In the present work, we continued the investigations of 
Ball et al.~\cite{ball2011,ball2015} about tracing delays, taking into account forward- and backward tracing. 
This study particularly focused on the questions 
how randomness in the tracing delay affects the effect 
of contact tracing to fight an epidemic during the initial phase of an outbreak, if forward- or backward tracing 
plays a more decisive role, and 
how the interplay between the time scales involved 
influences contact tracing. With respect to the last question, we focused on 
the mean tracing delay on the one hand, and the average 
time of infection (without control measure) respectively 
a latent period.\par
In order to approach these questions, we focused on the onset of 
an epidemic, and used a branching process approximation for the
spread of infections. On top of this linear branching process, 
contact tracing has been introduced. The tracing leads to
dependencies between individuals, which makes the stochastic 
process more complex
to analyze analytically. The main tool to analyze the process is the probability for an infected individual to be still infectious 
at a given age of infection. It has been possible to derive a system of integro-differential equations for this probability. As an explicit solution seems to be difficult to obtain, approximate solutions (for small tracing probability) have been derived. Based upon these approximations, an approximation for the reproduction number has been
proposed.
\par\medskip

The present study addressed the effect of randomness  in the tracing
delay: 
If we compare a fixed and an exponentially distributed tracing delay, we find 
that the effect of the fixed delay is smaller than that 
of the exponential distribution (if both have the same expectation). 
Most likely, the main reason for this observation is the exponential decrease of the probability to be infected 
at a certain age of infection: The expected number of infections produced after a certain age of infection will decrease exponentially (note that we do not condition on the fact
that a person reaches this age of infection). 
In a situation with randomness in the delay, 
some persons are detected earlier, and some later in comparison with the fixed delay. Due to the exponential
decrease, the gain in the effect by the more 
early detections 
is higher than the loss in effect by the later 
detections. This difference leads to 
an exponential 
respectively a polynomial decrease of the effect 
on the mean tracing delay. This finding may 
underline the importance 
to avoid outliers in the tracing delay: if the time between 
detection of an index case and investigation of contacts 
becomes large, contact tracing for this index case 
becomes inefficient. The efficiency of tracing 
for one index case 
does not decrease linearly in the tracing delay, but exponentially. 
Considering the high costs for contact tracing, 
it may be worth to implement a tracing program in such a way that 
a long tracing delay for an index cases is not likely to occur. 
However, a delay becomes notable only if it is in the 
range of the infectious period, as it is to expect. If it is 
distinctively shorter, the delay hardly plays a role, if it is longer, contact tracing becomes inefficient.\par\medskip

If we introduce a latency period, we find the same result for backward tracing (and a deterministic, fixed delay) as before: 
the effect decreases exponentially, where the ratio between tracing delay and latency period (delay in infectivity) is decisive. It is slightly different in forward tracing: here, 
the effect is even stronger compared with the case without 
any delay (no tracing delay and no latency period) if 
the tracing delay is shorter than the latency period. 
We clearly find a race between forward tracing and 
infections. Only if the tracing delay becomes larger 
than the latency period, we again find an exponential decrease in the effect. Forward tracing benefits from a long latency period, and may even become stronger than
backward tracing. If the tracing- and the latency period are in the same range, backward tracing is more likely to 
play the central role.\par\medskip

Basically, we have three ingredients of the implementation of 
a contact tracing program  that decide about its efficiency: 
(a) the probability for an infectious person to become an index case 
(b) how likely is a contact reported by an index case 
and (c) how large is the tracing delay. At lowest order, 
the effect (measured by the reproduction number) is the product of the 
first two probabilities times an exponentially decreasing function 
in the tracing delay. The time scale of this decrease is given by 
the recovery rate, and affected by the latency period. 
One may think about resource allocation within a tracing program. 
As long as
the tracing delay is distinctively shorter than the latency period
resp.\ the infectious period, it seems to be better to put effort 
in the detection of more contacts or index cases. Only if the 
disease is fast (short latency period/infectious period), it is 
of importance to decrease the tracing delay. If, however, the 
time scale of the infection is too fast, contact tracing as a control
measure could be inadequate. These findings are in line with, e.g., results
by Fraser et al.~\cite{fraser2004}

\bibliographystyle{abbrv} 
\bibliography{ctBib}

\begin{appendix}
\section{Latency period}
\label{latenecyAppendix}

We aim at a first order approximation (the only analysis that 
is feasible in general); we find
\begin{eqnarray*}
\kappa_0'(a) &=& 
-\kappa_0(a) \bigg(
\sigma(a) + \alpha(a)\nonumber\\
&&\qquad \quad
- p\,\int_0^a\phi(a-\tau)\int_0^\tau\bigg(\beta(c) 
       \widetilde\kappa'(c) +\beta(c) \alpha(c) \widetilde\kappa(c)\bigg)\, dc \, d\tau
\bigg)+{\cal O}(p^2)\\
\kappa_{i-1}(b)\,\kappa_i(a|b) & = & 
\kappa_{0}(a)\bigg\{
\kappa_{i-1}(b)\nonumber\\
&&\qquad
-p\, \int_0^a 
\left(
{- \widetilde\kappa}'(b+c)
-\,\alpha(b+c)\, \widetilde\kappa(b+c)\,
\right)\,\int_c^{a}\phi(a'-c)\, da'\, dc\bigg\}+{\cal O}(p^2)\\
\kappa_i(a) & = & \frac{\int_0^\infty 
\kappa_{i}(a|b)\beta(b)\kappa_{i-1}(b)\, db}
 {\int_0^\infty \beta(\tau)\kappa_{i-1}(\tau)\,d\tau }\\
 \widetilde\kappa(a) &=& e^{\int_0^a\,\sigma(\tau)+\alpha(\tau)\, d\tau}.
\end{eqnarray*}
 \par\medskip

\begin{prop} Up to second order in $p$, $\kappa_i(a)=\kappa_1(a)$ for $i>0$. 
\end{prop}
{\bf Proof:} We first find 
\begin{eqnarray}
&& \kappa_i(a|b) \nonumber\\
& = &
\kappa_{0}(a)\bigg\{
1 -\frac{p}{\kappa_{i-1}(b)}\, \int_0^a 
\left(
{- \widetilde\kappa}'(b+c)
-\,\alpha(b+c)\, \widetilde\kappa(b+c)\,
\right)\,\int_c^{a}\phi(a'-c)\, da'\, dc\bigg\}+{\cal O}(p^2)\nonumber\\
 & = & 
\kappa_{0}(a)\bigg\{
1
-p\, \int_0^a \,\,\,
\left(
\frac{{- \widetilde\kappa}'(b+c)}{\widetilde \kappa(b)}
-\,\alpha(b+c)\, \frac{\widetilde\kappa(b+c)}{\widetilde \kappa(b)}\,
\right)\,\int_c^{a}\phi(a'-c)\, da'\,\,\, dc\bigg\}+{\cal O}(p^2)
\label{xkappaab}
\end{eqnarray}
and hence $\kappa_i(a|b) = \kappa_1(a|b)+{\cal O}(p^2)$. 
We drop the index and simply write $\kappa(a|b) = \kappa_1(a|b)+{\cal O}(p^2)$, where $\kappa(a|b)$ is defined by 
equation (\ref{xkappaab}), where the ${\cal O}(p^2)$ terms are neglected. In zero order, 
$\kappa_0(a) = \widetilde\kappa(a)$, and hence also 
$\kappa(a|b)=\widetilde \kappa(a)$; 
we write 
$$\kappa(a|b) = \widetilde\kappa(a) + p \zeta(a|b)$$
where $\zeta(a|b)$ does not depend on $p$. Hence, 
\begin{eqnarray*}
\kappa_i(a) & = & \frac{\int_0^\infty 
\big(\widetilde\kappa(a)+p\zeta(a|b)\big)\beta(b)\kappa_{i-1}(b)\, db}
 {\int_0^\infty \beta(\tau)\kappa_{i-1}(\tau)\,d\tau }+{\cal O}(p^2)
 = \widetilde\kappa(a)
 + p \, \frac{\int_0^\infty 
\zeta(a|b)\beta(b)\kappa_{i-1}(b)\, db}
 {\int_0^\infty \beta(\tau)\kappa_{i-1}(\tau)\,d\tau }+{\cal O}(p^2)\\
  &=& \widetilde\kappa(a)
 + p \, \frac{\int_0^\infty 
\zeta(a|b)\beta(b)\widetilde\kappa(b)\, db}
 {\int_0^\infty \beta(\tau)\widetilde\kappa(\tau)\,d\tau }+{\cal O}(p^2).
 \end{eqnarray*}
 \par\qed\par\medskip

As a consequence, the 
first order correction of the asymptotic reproduction 
 number ($i\rightarrow\infty$) consist of two clearly 
 separated parts: that due to forward tracing, and that due 
 to backward tracing. Note that we choose the signs in the
 definition of $\eta_\pm$ below in such a way, that the 
 correction terms of the reproduction number 
 have a minus sign. 
 
 \begin{prop}
 Let $\kappa_0(a) = \widetilde\kappa(a)-p\eta_-(a)+{\cal O}(p^2)$, 
 $$
 \eta_+(a) = 
 \widetilde\kappa(a)
 \frac{\int_0^\infty 
 \int_0^a \,\,\,
\left(
\frac{{- \widetilde\kappa}'(b+c)}{\widetilde \kappa(b)}
-\,\alpha(b+c)\, \frac{\widetilde\kappa(b+c)}{\widetilde \kappa(b)}\,
\right)\,\int_c^{a}\phi(a'-c)\, da'\, dc\,\,
\beta(b)\widetilde\kappa(b)\, db}
 {\int_0^\infty \beta(\tau)\widetilde\kappa(\tau)\,d\tau }
 $$
 $r_\pm= \int_0^\infty \beta(a)\, \eta_\pm(a)\, da$. Then, 
 $$ R_{ct} = R_0 - p(r_++r_-)+{\cal O}(p^2).$$
 \end{prop}
 This proposition is a direct consequence of 
  equ.~(\ref{xkappaab}).
 
 \par\medskip
 
 Now we specify the parameters as indicated above:
$\beta(a)$, $\alpha(a)$ and $\sigma(a)=0$ for $a<t_i$, 
and are constant ($\beta$, $\alpha$, $\sigma$) afterwards; 
The contact tracing delay is a constant $T$. 
  We 
 compute the reproduction number and the first order correction terms for~$R_{ct}$. We will use
 $$\widehat\kappa(a) = e^{-(\alpha+\sigma)\,a},\quad
 \widetilde\kappa(a) = e^{-\int_0^a(\alpha(\tau)+\sigma(\tau))\,d\tau}
 $$
 such that $\widetilde\kappa(a)=\widehat\kappa(a-T_i)$ for $a>T_i$. 
 The zero order term is a direct
 consequence of $R_0 = \int_0^\infty \beta(a)\, \widetilde\kappa(a)\, da$.
 
\begin{prop}
$
R_0 = \frac{\beta}{\alpha+\sigma}$.
\end{prop}

\subsection{First order approximation: backward tracing} 
\label{latenecyAppendixBack}
 \begin{prop}
  The first order approximation reads
$$\kappa_0(a) = \widehat\kappa(a)+{\cal O}(p^2) \quad\mbox{for} \quad a\leq T_i+T_t.$$
and for $a>T_i+T$
\begin{eqnarray*}
\kappa_0(a)
&=& 
\widehat\kappa(a)\,
\,
\left[1-  p\,\beta\,p_{\text{obs}}\,
 \left( (a-T_i-T)-\frac1 {\alpha+\sigma}(1-\widetilde\kappa(a-T-T_i)\right)
 \right]
+ {\cal O}(p^2)
\end{eqnarray*}
 \end{prop}
{\bf Proof: } The proof consist of direct computations. 
\begin{eqnarray*}
\kappa_-'(a)
& = &
-\kappa_-(a) \left(
\sigma(a) + \alpha(a)
- p\,\int_0^a\phi(a-\tau)\int_0^\tau\bigg(\beta(\tau-c) 
       \widetilde\kappa'(c) +\beta(\tau-c) \alpha(c)\widetilde\kappa(c)\bigg)\, dc \, d\tau
\right)+ {\cal O}(p^2)\\
& = &
-\kappa_-(a) \left(
\sigma(a) + \alpha(a)
- p\,\int_0^a\phi(a-\tau)\int_0^{\tau-T_i}\bigg(\beta \sigma(c)\widetilde\kappa(c)\bigg)\, dc \,\chi_{\tau>T_i}\, d\tau
\right)+ {\cal O}(p^2)\\
& = &
-\kappa_-(a) \left(
(\sigma + \alpha)\,\,\chi_{a>T_i}
+ p\,\beta\,\sigma\,\int_0^a\phi(a-\tau)\int_{T_i}^{\tau-T_i}
       e^{(\alpha+\sigma)(c-T_i)}\, dc \,\,\,\chi_{\tau>2\,T_i}\, d\tau
\right)+ {\cal O}(p^2)\\
& = &
-\kappa_-(a) \left(
\sigma(a) + \alpha(a)
+ p\,\beta\,p_{\text{obs}}\,\int_0^a\phi(a-\tau)\,
       \left(1-\widehat\kappa(\tau-2\,T_i)\right)\,\,\chi_{\tau>2T_i}\, d\tau
\right)+ {\cal O}(p^2)\\
& = &
-\kappa_-(a) \bigg(
\sigma(a) + \alpha(a)
+ p\,\beta\,p_{\text{obs}}\,
       \left(1-\widehat\kappa(a-T-2\,T_i)\right)\,\,\chi_{a>2\,T_i+T}\bigg)\,
+ {\cal O}(p^2)
\end{eqnarray*}
where we used the definition $\widehat\kappa(a)=e^{-(\alpha+\sigma)a}$ for $a>0$. Then, 
\begin{eqnarray*}
\kappa_-(a) &=& \widehat\kappa(a)\,
\,
e^{\chi_{a>2\,T_i+T}\,  p\,\beta\,p_{\text{obs}}\,
 \left( (a-2\,T_i-T)-\frac 1 {\alpha+\sigma}(1-\widehat\kappa_-(a-T-T_i))\right)
}+ {\cal O}(p^2)\\
&=&
\widetilde\kappa(a)\,
\,
\left[1-  p\,\beta\,p_{\text{obs}}\,
\chi_{a>2\,T_i+T}\, 
 \left( (a-2\,T_i-T)-\frac 1 {\alpha+\sigma}(1-\widetilde\kappa_-(a-T-2\,T_i))\right)
 \right]
+ {\cal O}(p^2)
\end{eqnarray*}

\par\qed\par\medskip

The term $2\, T_i$ is a consequence of 
backward tracing: An infectee is only produced after $T_i$ 
time units, and the earliest time point at which this infectee 
can be observed is $T_i$ time units after his/her infection. Hence, 
the earliest time point at which a backward tracing event 
can be triggered is $2T_i$, and the infector is then traced 
at age of infection $2\,T_i+T$.\par\medskip

Consequently, we have (note our sign convention)
$$\eta_-(a) =  \widetilde\kappa(a)\,
\,
\,\beta\,p_{\text{obs}}\,
\chi_{a>2\,T_i+T}\, 
 \left( (a-2\,T_i-T)-\frac 1 {\alpha+\sigma}(1-\widehat\kappa(a-T-2\,T_i))\right).
$$
Evaluating the integral $r_-=\int_0^\infty\beta(a)\, \eta_-(a)\, da$ yields
\begin{eqnarray}
 r_- = \frac 1 2 \,\, p_{\text{obs}} \,\,R_0^2 \,\,\,\widehat\kappa(T+T_i).\end{eqnarray}

\subsection{First order approximation: forward tracing} 
\label{latenecyAppendixForw}
\begin{prop}
$\eta_+(a) = p_{\text{obs}}\,\,
\chi_{a>T}\,\,
\widetilde\kappa(a)(1-\widehat\kappa(a-T)).$
\end{prop}
{\bf Proof:}
We evaluate 
 $$
 \eta_+(a) = 
 \widetilde\kappa(a)
 \frac{\int_0^\infty 
 \int_0^a \,\,\,
\left(
\frac{{- \widetilde\kappa}'(b+c)}{\widetilde \kappa(b)}
-\,\alpha(b+c)\, \frac{\widetilde\kappa(b+c)}{\widetilde \kappa(b)}\,
\right)\,\int_c^{a}\phi(a'-c)\, da'\, dc\,\,
\beta(b)\widetilde\kappa(b)\, db}
 {\int_0^\infty \beta(\tau)\widetilde\kappa(\tau)\,d\tau }
 $$
 which yields the result for our special
 case as
 $ \int_0^\infty \beta(\tau)\widetilde\kappa(\tau)\,d\tau  = R_0,$ 
 $\int_c^a\phi(a'-c)\, da'=\chi_{a>T}\,\,\chi_{c<a-T}$, and
 $$
 \left(
\frac{{- \widetilde\kappa}'(b+c)}{\widetilde \kappa(b)}
-\,\alpha(b+c)\, \frac{\widetilde\kappa(b+c)}{\widetilde \kappa(b)}\,
\right)
=
\sigma(b+c)\, \frac{\widetilde\kappa(b+c)}{\widetilde \kappa(b)}.
$$
Note that the out-most integral only extends 
over $b\in[T_i,\infty)$ as $\beta(b)=0$ for $b<T_i$.  
For $b>T_i$, we have 
$$\sigma(b+c)\,\widetilde\kappa(b+c)/\widetilde \kappa(b) 
= \sigma\,\widehat\kappa(b+c-T_i)/\widehat \kappa(b-T_i) 
= \sigma\,\widehat \kappa(c).$$
That is, the inner integral does not depend on $b$, and the result follows easily.  
\par\qed\par\medskip

 Integrating $\beta(a)\, \eta_+(a)$ over $a\in\R_+$ 
 yields the correction term for forward tracing, 
 \begin{eqnarray}
r_+ &=& \frac 1 2\,\,
p_{\text{obs}} R_0 \,\,
\widehat\kappa(\max(T_i, T)-T_i)\,\,
\bigg(
2
-\,\,
\widehat\kappa(\max(T_i, T)-T)
\bigg)\nonumber\\
 &=& \frac 1 2\,\,
p_{\text{obs}} R_0 \,\,
\frac{\widehat\kappa(\max(T_i, T))}
{\widehat\kappa(T_i)}\,\,
\bigg(
2
-\,\,
\frac{\widehat\kappa(\max(T_i, T))}
{\widehat\kappa(T)}
\bigg).
 \end{eqnarray}

\end{appendix}
\end{document}